\def\ove#1{\overline{#1}}    %
\def\ovs#1#2{\overset{#1}{#2}}        

     \def\({\left(}                  \def\lee{\leqslant}
     \def\){\right)}      \def\e{\varepsilon}    
     \def\[{\left[}       
     \def\]{\right]}      \def\ffi{\varphi}

                                      \def\ot{\otimes}
     \def\<{\langle}                 \def\wh{\widehat}
     \def\>{\rangle}                 \def\wt{\widetilde}
                 \def\sbs{\subset}
\def\tr{\operatorname{trace}\,}

\def\N{\operatorname{N}}
\def\I{\operatorname{I}}

\def\QN{\operatorname{QN}}

\def\R{\operatorname{R}}

                     \def\sbs{\subset}

\documentclass[12pt,oneside]{amsproc}
\usepackage[T2A]{fontenc}
\usepackage[english]{babel}
\usepackage{amssymb}
\usepackage{amsthm}

\title{Banach spaces without approximation properties of type $p$}
\author{Oleg Reinov${ }^*$\, and Qaisar Latif}
\address{Oleg Reinov:\newline Department of Mathematics and Mechanics, St. Petersburg State University,
Saint Petersburg, RUSSIA.\newline
\phantom{Ao} Abdus Salam School of Mathematical Sciences, 68-B, New Muslim Town, Lahore 54600, PAKISTAN.
}
\email{orein51@mail.ru}
\address{Qaisar Latif:\newline Abdus Salam School of Mathematical Sciences, 68-B, New Muslim Town, Lahore 54600, PAKISTAN.}

\thanks{${ }^*$The research was supported by the Higher  Education Commission of Pakistan.}

\thanks{${ }$ AMS Subject Classification 2000: 46B28.
Spaces of operators; tensor products; approximation properties }
\thanks{${ }$ Key words: $p$-nuclear operators, approximation
properties, compact approximation, nuclear tensor norms, tensor
products. }

\begin{document}
\vphantom{} \maketitle

\begin{abstract}
The main purpose of this note is to show that the question posed in the paper
of Sinha D.P. and Karn A.K.
("Compact operators which factor through subspaces of $l_p$
Math. Nachr.  281, 2008, 412-423; see the very end of that paper)
has a negative answer, and that the answer could be obtained,
essentially, in 1985 after the papers [2], [3] by Reinov O.I.
have been appeared in 1982 and in 1985 respectively.
\end{abstract}

\vskip 0.3cm

The main purpose of this note is to show that the question posed in [1]  
(see the very end of that paper) has a negative answer, and that the answer could be obtained,
 essentially, in 1985 after the papers [2], [3] have been appeared 
in 1982 and in 1985 respectively.

For the sake of completeness (and since the article [3] is now quite   
difficult to get)  
we will reconstruct some of the facts from [3], 
maybe, in a somewhat changed form (a translation of [3] with some remarks
can be found in  arXiv:1002.3902v1 [math.FA]).             

All the notation and terminology, we use, are from [2], [4], or [5]
and more or less standard. 
For the general theory of absolutely $p$-summing, $p$-nuclear and
other operator ideals, we refer to [6].           
We will keep, in particular, the following standard notations.
If $A$ is a bounded subset of a Banach space $X$, then
$\ove{\Gamma(A)}$ is the closed absolutely convex hull of $A;$
$ X_A$ is the Banach space with "the unit ball" $ \ove{ \Gamma(A)}$;
$ \Phi_A: X_A\to X$ is the canonical embedding. For $ B\sbs X,$
it is denoted by $ \overline{B}^{\,\tau}$ and $ \ove{B}^{\,\|\cdot\|}$
the closures of the set $ B$ in the topology $ \tau$ and in the norm $ \|\cdot\|$
respectively. When it is necessary, we denote by $ \|\cdot\|_X$
the norm in $X.$ Other notations:\,
$ \Pi_p,$ $ \QN_p,$ $ \N_p,$ $ \I_p$ are the ideals of absolutely $p$-summing,
quasi-$p$-nuclear, $p$-nuclear, strictly $p$-integral operators,
respectively;
$ X^*\wh\ot_p Y$ is the complete tensor product associated with
$ \N_p(X,Y);$
$ X^*\wh{\wh\ot}_p Y= \ove{X^* \ot Y}^{\,\pi_p}$ \, (i.e. the closure of the set
of finite dimensional operators in $ \Pi_p(X,Y)).$ Finally,
if $ p\in[1,+\infty]$ then $ p'$ is the adjoint exponent.

Now, we recall the main definition from [3]  
of the topology $\tau_p$ of the $\pi_p$-compact convergence, which is
just the same as the $\lambda_p$-topology in [1]\footnote{         
  Definition 4.3 from [1]:\,   
 "Given a compact subset $K\sbs X$ we define a seminorm $||\cdot||_K$
on $\Pi_p(X,Y)$ given by
$||T||_K=\inf\{\kappa_p^d(Ti_Z):\, i_Z:Z\to X \text{ as above}\}.$
The family of seminorms $\{||\cdot||_K:\, K\sbs X \text{ is compact}\}$\,
determines a locally convex topology $\lambda_p$ on $\Pi_p(X,Y).$"; ---
for details, and for the definition of maps $i_Z,$ see [1], Lemma 4.2.
  },
where it was proved, essentially, that that
$T\in K_p^d(X,Y)$ if and only if $T\in QN_p(X,Y).$

For Banach spaces $X, Y,$ the {\it topology $ \tau_p$\ of
$ \pi_p$\!-compact convergence}\
in the space $ \Pi_p(Y,X)$ is the topology, a local base (in zero) of which
is defined by sets of type
$$ \omega_{K,\e}= \left\{ U\in \Pi_p(Y,X):\ \pi_p(U\Phi_K)<\e\right\},
$$
where $ \e>0,$\, $ K=\ove{\Gamma(K)}$ --- a compact subset of $ Y.$

\vskip 0.2cm
{\bf Proposition 1} [3].\ {\it
Let $ \operatorname{R}$ be a linear subspace in
$ \Pi_p(Y,X),$  containing
$ Y^*\ot X.$ Then
$ (\operatorname{R},\tau_p)'$ is isomorphic to a factor space of the space
$ X^*\wh\ot_{p'} Y.$
More precisely, if $ \ffi\in (\operatorname{R}, \tau_p)',$
then there exists an element $ z=\sum_1^\infty x'_n\ot y_n\in X^*\wh\ot_{p'} Y$
such that
$$ \ffi(U)= \tr\, U\circ z,\ \, U\in \operatorname{R}. \eqno(*)
$$
On the other hand, for every
$ z\in X^*\wh\ot_{p'} Y$ the relation
$ (*)$ defines a linear continuous functional on
$ (\R,\tau_p).$
  }
\vskip 0.1cm

{\it Proof.}\  
Let
$ \ffi$ be a linear continuous functional on $ (\R,\tau_p).$
Then one can find a neighborhood of zero
$ \omega_{K}=\omega_{K,\e},$ such that $ \ffi$ is bounded on it:
$ \forall\, U\in\omega_K, \ |\ffi(U)|\lee 1.$
We may assume that $ \e=1.$ Consider the operator
$ U\Phi_K:\,Y_K\ovs{\Phi_K}\longrightarrow  Y
               \ovs{U}\longrightarrow X.$
Since the mapping  $ \Phi_K$ is compact,  $ U\Phi_K\in \QN_p(Y_K,X).$
Put $ \ffi_K(U\Phi_K)=\ffi(U)$ for $ U\in \R.$
On the linear subspace
$ \R_K= \left\{ V\in \QN_p(Y_K,X):\ V=U\Phi_K \right\}$
of the space  $ \QN_p(Y_K,X),$
the linear functional
$ \ffi_K$
is bounded: if
$ V=U\Phi_K\in \R_K$ and $ \pi_p(V)\lee 1,$
then $ |\ffi_K(V)|=|\ffi(U)|\lee 1.$
Therefore, $ \ffi_K$ can be extended to a linear continuous functional
$ \wt\ffi$ on the whole  $ \QN_p(Y_K,X);$ moreover, because of the injectivity of the ideal
$ \QN_p,$
considering $ X$ as a subspace of some space $ C(K_0),$
we may assume that $ \wt\ffi\in\QN_p(Y_K, C(K_0))^*.$ Let us mention that
$$ \wt\ffi(jU\Phi_K)=\ffi_K(U\Phi_K)=\ffi(U) \eqno(1)
$$
(here $ j$ is an isometric embedding  of $ X$ into $ C(K_0)$).

Furthermore, since
$ \QN_p(Y_K, C(K_0))^*= \I_{p'}(C(K_0), (Y_K)^{**}),$
we can find an operator
$ \Psi: C(K_0)\to (Y_K)^{**},$
for which
$$ \wt\ffi(A)= \tr \Psi A,\ \, A\in (Y_K)^*\ot C(K_0).
$$

Let
$ A_n\in (Y_K)^*\ot C(K_0),\, $ $ \pi_p(A_n-jU\Phi_K)\to 0.$
Then
$$ \wt\ffi(jU\Phi_K)= \lim\, \tr \Psi A_n. \eqno(2)
$$

Consider the operator
$ \Phi_K^{**}\Psi: C(K_0)\ovs{\Psi}\longrightarrow (Y_K)^{**}
   \ovs{\Phi_K^{**}}\longrightarrow Y.$
Since
$ \Psi\in \I_{p'},$ and $ \Phi_K$
is compact, we have
$ \Phi_K^{**}\Psi \in \N_{p'}(C(K_0), Y)= C(K_0)^*\wh\ot_{p'} Y.$
Let
$ \sum_1^\infty \mu_n\ot y_n\in C(K_0)^*\wh\ot_{p'} Y$
be a representation of the operator
$ \Phi^{**}_K\Psi.$
Put
$ z=\sum j^*(\mu_n)\ot y_n.$
The element
$ z$
generates an operator
$ \Phi_K^{**}\Psi j$
from
$ X$
to
$ Y.$
We will show now that
$ \tr U\circ z=\wt\ffi(jU\Phi_K)$
(note that
$ U\circ z$
is an element of the space
$ X^*\wh\ot X,$
so the trace is well defined).
We have:
$$  \tr U\circ z= \tr \( \sum j^*(\mu_n)\ot Uy_n\)=
   \sum \< j^*(\mu_n), Uy_n\> = $$
$$\sum \< \mu_n, jUy_n\>=
\tr jU\Phi_K^{**}\Psi=\tr (jU\Phi_K)^{**}\Psi, \eqno(3)
$$
where
$ (jU\Phi_K)^{**}\Psi:
    C(K_0)\ovs{\Psi}\longrightarrow (Y_K)^{**}
         \ovs{\Phi_K^{**}}\longrightarrow Y
       \ovs{U}\longrightarrow X\ovs{j}\longrightarrow C(K_0).$
Since
$ \pi_p(A_n- jU\Phi_K)\to 0,$
then
$ \pi_p\(A^{**}_n - (jU\Phi_K)^{**}\)\to 0.$
Moreover, if
$A:= A_n=\sum_1^N w_m\ot f_m\in (Y_K)^*\ot C(K_0),$
then
$$ \tr A_n^{**}\Psi= \sum_m \< \Psi^*w_m, f_m\> =
   \sum_m \< w_m, \Psi f_m\>= \tr \Psi A.
$$

Hence,
$ \tr (jU\Phi_K)^{**}\Psi = \lim\,\tr A^{**}_n\Psi=
\lim\,\tr \Psi A_n.$
Now, it follows from (3) and (2) that
$ \wt\ffi(jU\Phi_K)=\tr U\circ z.$
Finally, we get from (1):
$ \ffi(U)= \tr U\circ z.$
Thus, the functional
$ \ffi$
is defined by an element of
$ X^*\wh\ot_{p'} Y.$

Inversely, if
$ z\in X^*\wh\ot_{p'} Y,$
put
$ \ffi(U)=\tr U\circ z$
for
$ U\in \R$
(the trace is defined since
$ U\circ z\in X^*\wh\ot X).$
We have to show that the linear functional
$ \ffi$
is bounded on a neighborhood
$ \omega_{K,\e}$
of zero in
$ \tau_p.$
For this, we need the following fact,
which proof is rather standard  (see [3], Lemma 1.2):   
{\it
If
$ z\in X^*\wh\ot_q Y,$
then $ z\in X^*\wh\ot_q Y_K,$
where
$ K=\ove{\Gamma(K)}$
is a compact in
$ Y.$
}\
Now, let
$ K$ be a compact subset of
$ Y,$
for which
$ z\in X^*\wh\ot_{p'} Y_K.$
If
$ U\in \omega_{K,1},$
then
$ \pi_p(U\Phi_K)<1$
and
$ |\tr U\Phi_K\circ z|\lee
\|z\|_{X^*\wh\ot_{p'} Y_K}\cdot$  
$\pi_p(U\Phi_K)\lee C.$


\vskip 0.1cm

{\bf Corollary 1.}\it \, $ (\R,\tau_p)'=(\R,\sigma)',$
where
$ \sigma=\sigma(\R, X^*\wh\ot_{p'} Y).$
Thus, the closures
of convex subsets of the space
$ \Pi_p(Y, X)$
in
$ \tau_p$
and in
$ \sigma$
are the same.
\rm

\vskip 0.1cm

{\bf Proposition 2} [3]. \it
If the canonical mapping
$ j: X^*\wh\ot_{p'} Y\to \N_{p'}(X,Y)$
is one-to-one then
$ \Pi_p(Y,X)= \ove{Y^*\ot X}^{\,\tau_p}.$
\vskip 0.1cm

\rm
{\it Proof}
If the map
$ j$
is one-to-one then the annihilator
$ j^{-1}(0)^\perp$
of its kernel in the space, dual to
$ X^*\wh\ot_{p'}Y,$
coincides with
$ \Pi_p(Y, X^{**}).$
On the other hand, in any case
$ j^{-1}(0)^\perp=
\ove{Y^*\ot X}^{\,*}$
(the closure in
$ {}^*$\!-weak topology of the space
$ \Pi_p(Y, X^{**}));$
by Corollary from Proposition 1,
$$ \Pi_p(Y,X)\cap \ove{Y^*\ot X}^{\,*}= \ove{Y^*\ot X}^{\,\tau_p}.
$$

Therefore,
$ \Pi_p(Y,X)= \ove{Y^*\ot X}^{\,\tau_p}.$
\vskip 0.1cm

For a reflexive space
$ X,$
the dual space to
$ X^*\wh\ot_{p'} Y$
is equal to
$ \Pi_p(Y,X).$
Consequently, it follows from the last two statements

\vskip 0.1cm

{\bf Corollary 2.}\it \, For a reflexive space
$ X$
the canonical mapping
$ j: X^*\wh\ot_{p'} Y\to \N_{p'}(X,Y)$
is one-to-one iff the set of finite rank operators is dense in the space
$ \Pi_p(Y,X)$
in the topology
$ \tau_p$\ of $ \pi_p$-compact convergence. \rm

\vskip 0.1cm

Recall a definition of the {\it approximation property $AP_p$ of order}\, $p, p\in (0,\infty]$
(see, e.g., [2] and [7]):    
a Banach space $X$ has the $AP_p$ if for every Banach space $Y$
(equivalently, for every reflexive Banach space $Y,$ see [2] or [5])
one has the equality
 $Y^*\wh\ot_p X = \N_p(Y,X).$ It follows now from Proposition 1 and Corollary 2:

\vskip 0.1cm
{\bf Corollary 3.}\it \, For $p\in[1,\infty]$ and for every Banach space
$ X$
the following are equivalent:

$1)$ $X$ has the $AP_p;$

$2)$
for each Banach space
$ Y$\
$ \ove{X^*\ot_{p'} Y}^{\,\tau_{p'}}= \Pi_{p'}(X,Y);$

$3)$
for each reflexive Banach space
$ Y$\
$ \ove{X^*\ot_{p'} Y}^{\,\tau_{p'}}= \Pi_{p'}(X,Y);$

$3')$
for each reflexive Banach space
$ Y$\
$ \ove{X^*\ot_{p'} Y}^{\,\lambda_{p'}}= \Pi_{p'}(X,Y).$

\vskip 0.1cm
\rm

As shown in [2], [4] and [5] , for each $p, p\in [1,\infty], p\neq2,$
there exist reflexive Banach spaces without the $AP_p.$
Thus, we get from the last consequence

\vskip 0.1cm
{\bf Corollary 4.}\it \, For $1\le p\le \infty, p\neq2,$
there are two reflexive Banach spaces $X$ and $Y$ such that the natural map
$Y^*\wh\ot_p X \to \N_p(Y,X)$ is not one-to-one and $ \ove{X^*\ot_{p'} Y}^{\,\tau_{p'}}\neq \Pi_{p'}(X,Y).$
\vskip 0.1cm
\rm

Since, as we said, $\tau_q$ equals to
the topology $\lambda_q$ from [1], we get firstly Theorem 4.11 [1] (the case $p>2$ below)
and secondly the negative answer to
the question in the end of the paper [1] (whether for $p\in[1,2)$ 
every Banach space has the approximation property
of type $p$ of [1]\footnote{
By [1], a Banach space $X$ is said to have the {\it approximation property
of type $p$}\, if for every Banach space $Y,$ the finite rank operators
$\mathcal F(Y,X)$ is dense in $\Pi_p(Y,X)$ in the $\lambda_p$-topology.}). 
Or, generally, we have:

\vskip 0.1cm
 {\bf Theorem.}\it \, Let $1\le p\neq2 \le\infty.$ Then there is a $($reflexive$)$ Banach space
that fails to have the approximation property of type $p$
of $[1].$

\rm


\bigskip


\begin{thebibliography}{0}

\medskip



\bibitem{1} Sinha D.P., Karn A.K.:
\textit{Compact operators which factor through subspaces of $l_p$},
Math. Nachr.  281(2008), 412-423.

\bibitem{2} O.I. Reinov:
\textit{Approximation properties of order p and the existence of
 non-p-nuclear operators with p-nuclear second adjoints},
Math. Nachr.  109(1982), 125-134.



\bibitem{3} O.I. Reinov:
\textit{Approximation of operators in Banach spaces},
Application of functional analysis in the approximation theory (KGU, Kalinin)
  (1985), 128-142.

\bibitem{4} O.I. Reinov:
\textit{Disappearance of tensor elements in the scale of p-nuclear operators},
Theory of operators and theory of functions (LGU) 1(1983), 145-165.


\bibitem{5} O. Reinov: \textit{Approximation properties of order p and the existence of
  non-p-nuclear operators with p-nuclear second adjoints},
  Doklady AN SSSR, 256(1981), 43--47.



\bibitem{6} A. Pietsch:
Operator ideals, North-Holland, Deutscher Verlag der Wiss., Berlin,
1978.

\bibitem{7} P. Saphar:
\textit{Produits tensoriels d'espaces de Banach et classes
d'applications lineaires},
Studia Math. 38(1970), 71--100.



\end{thebibliography}
\end{document}